%% file: logdiff-preprint.tex
\documentclass{article}
\usepackage{amsmath,amsfonts,amsthm,mathrsfs,amscd,amssymb,pb-diagram}
\usepackage[letterpaper]{geometry}

\title{Logarithmic Combinatorial Differentials}
\author{Daniel Schepler}

\include{common-defs}

\begin{document}
\maketitle

\begin{abstract}
  Given a morphism $X \to S$ of fine log schemes, we develop a geometric
  description of the sheaves of higher-order differentials $\Omega^n_{X/S}$
  for $n > 1$, as well as a definition of the de~Rham complex in terms
  of this description.
\end{abstract}

\tableofcontents

\include{intro}
\include{diff}
\include{derham}
\include{coeff}
\section*{Acknowledgements}

The author wishes to thank Lawrence Breen for inviting him to
l'Universit\'e Paris 13 to give a talk on preliminary work on this
subject, for useful discussions during that visit, and for providing
comments on a draft version of this article.

The author has been supported by the research network {\em Arithmetic
  Algebraic Geometry} of the European Community (Contract
MRTN-CT-2003-504917).
\bibliographystyle{amsalpha}
\bibliography{logdiff}
\end{document}

%% file: common-defs.tex
\newcommand\bbZ{\mathbb{Z}}
\newcommand\scrO{\mathscr{O}}
\newcommand\gp{\mathop{\operatorname{gp}}}
\newcommand\Spec{\mathop{\operatorname{Spec}}}
\newcommand\dlog{\mathop{d\log}}
\newcommand\wedgetilde{\mathbin{\tilde\land}}
\newcommand\id{\mathop{\operatorname{id}}}
\newcommand\scrM{\mathscr{M}}
\newcommand\MX{\scrM_X}
\newcommand\MXgp{\scrM_X^{\gp}}
\newtheorem{theorem}{Theorem}[section]
\newtheorem{lemma}[theorem]{Lemma}
\newtheorem{proposition}[theorem]{Proposition}
\newtheorem{corollary}[theorem]{Corollary}
\theoremstyle{remark}
\newtheorem{remark}[theorem]{Remark}


%% file: intro.tex
\section*{Introduction}

Given a smooth morphism $X \to S$ of schemes, it is standard to define
$\Omega^1_{X/S} := I / I^2$, where $I$ is the ideal sheaf of the diagonal
in $X \times_S X$.  One normally then defines $\Omega^q_{X/S} := \bigwedge^q \Omega^1_{X/S}$
for $q > 1$.  On the other hand, in \cite{breen-messing}, Breen and
Messing give an alternate definition of $\Omega^q_{X/S}$ extending the
geometric definition of $\Omega^1_{X/S}$.  This paper was inspired by
similar definitions introduced by A.~Kock in his study of synthetic
differential geometry \cite{kock}, which in turn was an attempt to
transpose the methods in algebraic geometry, due to Grothendieck and
others, of studying the concept of infinitesimally close points to the
setting of $C^\infty$-manifolds.

For simplicity, let us assume that 2 is invertible on $S$.  Let $\Delta^n_{X/S}
:= X \times_S X \times_S \cdots \times_S X$ be the $n+1$-fold product, with the factors
indexed from 0 to $n$.  For $0 \leq i, j \leq n$, let $I_{ij}$ be ideal of
$\scrO_{\Delta^n_{X/S}}$ defining the partial diagonal $\{ (x_0, \ldots, x_n) \in \Delta^n_{X/S} :
x_i = x_j \}$.  Now let $\Delta^{(n)}_{X/S}$ denote the closed subscheme of $\Delta^n_{X/S}$
defined by $\sum_{0\leq i,j \leq n} I_{ij}^2$, and $\tilde I_{ij}$ the image of
$I_{ij}$ in $\scrO_{\Delta^{(n)}_{X/S}}$.  Then
\[ \prod_{i=1}^n \tilde I_{i-1,i} = \bigcap_{i=1}^n \tilde I_{i-1,i} = \prod_{i=1}^n
\tilde I_{0i} = \bigcap_{i=1}^n \tilde I_{0i} = \bigcap_{0 \leq i,j \leq n} \tilde
I_{ij}, \] and this common ideal, considered as an $\scrO_X$-module
via any of the $n+1$ projections $\Delta^{(n)}_{X/S} \to X$, is canonically
isomorphic to $\Omega^n_{X/S}$.  (In the general case, this construction
instead gives the $n$th antisymmetric power of $\Omega^1_{X/S}$.)

Our first observation is that in the general case, we can fix this
discrepancy by starting with the divided power envelope $D(n)$ of the
diagonal in $\Delta^n_{X/S}$.  In other words, if we let $\Delta^{[n]}_{X/S}$ be the closed
subscheme of $D(n)$ defined by $\sum_{0 \leq i,j \leq n} \bar I_{ij}^{[2]}$,
and $\tilde I_{ij}$ the image of $\bar I_{ij}$ in $\scrO_{\Delta^{[n]}_{X/S}}$,
then the five ideals above are once again equal, and are canonically
isomorphic to $\Omega^n_{X/S}$.  (In \cite{breen-messing}, Breen and Messing
corrected the discrepancy by expanding $\sum_{0 \leq i,j \leq n} I_{ij}^2$ in a
non-symmetric way.)

Log geometry provides a convenient language for discussing topics
related to compactification and singularities.  Recall that a pre-log
scheme $X$ is a scheme $X$ equipped with a sheaf of commutative
monoids $\MX$ and a morphism $\alpha_X : \MX \to \scrO_X^\times$, where
$\scrO_X^\times$ is the multiplicative monoid of $\scrO_X$.  (Note that we
use additive notation for $\MX$, thus considering $m \in \MX$ to be a
logarithm of $\alpha(m)$, and considering $\alpha$ to be an exponentation map.)
This is a log scheme if the induced morphism $\alpha_X^{-1}(\scrO_X^*) \to
\scrO_X^*$ is an isomorphism.  A log scheme is called {\em fine} if
locally the log structure is induced by a pre-log structure $P \to
\scrO_X^\times$ where $P$ is the constant sheaf of a finitely-generated
integral monoid.  Given a morphism $X \to S$ of log scheme, Kato
\cite{kato} defines a universal sheaf of relative log differentials
$\Omega^1_{X/S}$ with a log derivation $(d, \dlog) : (\scrO_X, \scrM_X) \to
\Omega^1_{X/S}$.  This means that $d : \scrO_X \to \Omega^1_{X/S}$ is an
$\scrO_S$-derivation, $\dlog : \MX \to \Omega^1_{X/S}$ is an additive map
annihilating the image of $\scrM_S$, and for $m \in \MX$, we have
\[ d \alpha(m) = \alpha(m) \dlog m. \]
For example, suppose $X$ is a smooth scheme over a field $k$, and $D$
is a divisor with normal crossings on $X$.  Let $Y := X \setminus D$, with
open immersion $i : Y \to X$.  We then define $\MX := i_* \scrO_Y^* \cap
\scrO_X$, with $\alpha_X$ the natural inclusion map.  This defines a log
scheme, and the sheaf of log differentials $\Omega^1_{X/k}$ is exactly the
classical sheaf $\Omega^1_{X/k}(\log D)$ of differentials with log poles
along $D$.

Our aim in this paper is to extend Breen and Messing's theory to give
an intrinsic geometric description of $\land^n \Omega^1_{X/S}$ for $n>1$ in the
case of log schemes.  Thus, consider a morphism $X \to S$ of fine log
schemes.  (Note that we do not require this morphism to be log
smooth.)  Again, let $\Delta^n_{X/S} := X \times_S \cdots \times_S X$ be the $n+1$-fold
product.  Then there exists a right universal log scheme $D(n)$ with
an {\em exact} closed immersion $X \to D(n)$ defined by a PD ideal on
$D(n)$, and a morphism $D(n) \to \Delta^n_{X/S}$, factoring the diagonal
morphism $X \to \Delta^n_{X/S}$ \cite{kato}.  Again, let $\Delta^{[n]}_{X/S}$ be
the closed subscheme of $D(n)$ defined by the ideal $\sum_{0 \leq i,j \leq n}
\bar I_{ij}^{[2]}$, where $\bar I_{ij}$ is the ideal of the partial
diagonal $\{ x_i = x_j \}$ in $\Delta(n)$, and $\tilde I_{ij}$ the image of
$\bar I_{ij}$ in $\scrO_{\Delta^{[n]}_{X/S}}$.  Then we will prove that in
this more general case, once again the five ideals above are equal and
are canonically isomorphic to $\Omega^n_{X/S}$.  The proof we give here is
an improvement on the proof given in \cite{breen-messing}.

In terms of this description, the de~Rham complex becomes particularly
simple, in the form of an Alexander-Spaniel complex.  First, for $m, n
\geq 0$, consider $\Delta^m_{X/S}$ as a scheme over $X$ via the last
projection, and $\Delta^n_{X/S}$ as a scheme over $X$ via the first
projection.  Then we have a morphism
\begin{align*}
\Delta^{m+n}_{X/S} & \to \Delta^m_{X/S} \times_X \Delta^n_{X/S}, \\
(x_0, \ldots, x_m, \ldots, x_{m+n}) & \mapsto ((x_0, \ldots, x_m),
(x_m, \ldots, x_{m+n})).
\end{align*}
This induces a map $\Delta^{[m+n]}_{X/S} \to \Delta^{[m]}_{X/S} \times_X \Delta^{[n]}_{X/S}$, which in turn induces the
wedge product.  Similarly, given $n \geq 0$ and $0 \leq i \leq n+1$, define $d_i :
\Delta^{n+1}_{X/S} \to \Delta^n_{X/S}$ to be the map which forgets the $i$th component.  This
induces maps $d_i : \Delta^{[n+1]}_{X/S} \to \Delta^{[n]}_{X/S}$, and the differential $d :
\Omega^n_{X/S} \to \Omega^{n+1}_{X/S}$ is induced by
\[ d_0^* - d_1^* + \cdots + (-1)^{n+1} d_{n+1}^* : \scrO_{\Delta^{[n]}_{X/S}} \to
\scrO_{\Delta^{[n+1]}_{X/S}}. \]

Finally, suppose $X \to S$ is log smooth.  We observe that each
$\Delta^{[n]}_{X/S}$ is an object in the log crystalline site of $X$ over
$S$, and each $d_i$ is a morphism in this site.  Therefore, given a
crystal $E$ on this site, which corresponds to a module with
quasi-nilpotent connection $(E_X, \nabla)$, we have transition maps
$\theta_{d_i} : E_{\Delta^{[n]}_{X/S}} \to E_{\Delta^{[n+1]}_{X/S}}$.  Here
$E_{\Delta^{[n]}_{X/S}} \simeq E_X \otimes_{\scrO_X} \scrO_{\Delta^{[n]}_{X/S}}$ via the
isomorphism $\theta_{\pi_0} : \pi_0^* E_X \to E_{\Delta^{[n]}_{X/S}}$.  We will show
that the differential \[ \nabla : E_X \otimes_{\scrO_X} \Omega^n_{X/S} \to E_X
\otimes_{\scrO_X} \Omega^{n+1}_{X/S} \] in the de~Rham complex of $(E_X, \nabla)$ is
induced by
\[ \theta_{d_0} - \theta_{d_1} + \cdots + (-1)^{n+1} \theta_{d_{n+1}} : E_X \otimes_{\scrO_X}
\scrO_{\Delta^{[n]}_{X/S}} \to E_X \otimes_{\scrO_X} \scrO_{\Delta^{[n+1]}_{X/S}}. \]

%% file: diff.tex
\section{Combinatorial Differentials}

\subsection{Local Construction}

We begin in this section with a simplified situation: suppose $A$ is a
ring, $B$ an $A$-algebra, and $Q \to P$ a morphism of finitely-generated
integral monoids with compatible maps $Q \to A$ and $P \to B$.  Let $S :=
\Spec A$, $X := \Spec B$, with the log structures induced by $Q$ and
$P$, respectively.  Now let $\Delta^n_{X/S} := X \times_S X \times_S \cdots \times_S X$ be the
$n+1$-fold product, with the factors indexed from 0 to $n$.  In other
words, $\Delta^n_{X/S} = \Spec (B_n)$, where $B_n := B \otimes_A \cdots \otimes_A B$, with log
structure induced by $P_n$, the quotient of $P \oplus P \oplus \cdots \oplus P$ by the
congruence generated by
\[ (q, 0, \cdots, 0) \equiv (0, q, \cdots, 0) \equiv \cdots \equiv (0, 0, \cdots, q) \]
for $q \in Q$.  Then $P_n^{\gp}$ is the quotient of $P^{\gp} \oplus P^{\gp} \oplus
\cdots \oplus P^{\gp}$ by $\{ (q_0, q_1, \ldots, q_n) \in Q^{\gp} \oplus \cdots \oplus Q^{\gp} : q_0 +
\cdots + q_n = 0 \}$.

Now the diagonal map $X \to \Delta^n_{X/S}$ corresponds to the product map $B_n \to
B$, and it has a chart given by the sum map $P_n \to P$.  Now
let
\[ P_n' := \{ (p_0, p_1, \ldots, p_n) \in P_n^{\gp} : p_0 + p_1 + \ldots + p_n \in P
\}, \] $B_n' := B_n \otimes_{\bbZ[P_n]} \bbZ[P_n']$, and $Z_n := \Spec B_n'$
with the log structure induced by $P_n'$.  (For $p \in P$, we will use
the notation $e^p$ for the corresponding element of $\bbZ[P]$, in
order to avoid confusion between addition in $P$ and addition in
$\bbZ[P]$.)  Then the map $X \to Z_n$ corresponding to the sum map $P_n'
\to P$ is an exact closed immersion, and the map $Z_n \to \Delta^n_{X/S}$
corresponding to the inclusion $P_n \hookrightarrow P_n'$ is log \'etale.  Therefore,
$Z_n$ may be used as the basis for constructing the log infinitesimal
neighborhoods and the divided power envelope of $X$ in $\Delta^n_{X/S}$
\cite{kato}.  (Recall that a map $g : Q \to P$ of integral monoids is
{\em exact} if $(g^{\gp})^{-1}(P) = Q$, and a morphism $X \to S$ of fine
log schemes is exact if for every point $x \in X$ with image $s \in S$,
$\scrM_{S,s} \to \scrM_{X,x}$ is exact.  A log closed immersion $f : X \to
Y$ is exact if and only if it is strict, i.e.~$f^* \scrM_Y \to \scrM_X$
is an isomorphism, where $f^* \scrM_Y$ is the log structure induced by
$f^{-1} \scrM_Y$.)




For notation, let $\pi_i^* : P \to P_n'$ be the $i$th inclusion map,
corresponding to the $i$th projection $\pi_i : Z_n \to X$.
Now for each pair $i,j$ with $0 \leq i,j \leq n$, we have a closed immersion
$m_{ij} : Z_{n-1} \to Z_n$ corresponding to the map $\mu_{ij} : B_n' \to B_{n-1}'$,
\begin{align*}
  (y_0 \otimes \cdots \otimes y_i \otimes \cdots \otimes y_j \otimes \cdots \otimes y_n) & \otimes e^{(p_0, \ldots, p_i, \ldots, p_j, \ldots,
  p_n)} \mapsto \\
  (y_0 \otimes \cdots \otimes y_i y_j \otimes \cdots \otimes \hat y_j \otimes \cdots \otimes y_n) & \otimes e^{(p_0, \ldots, p_i +
    p_j, \ldots, \hat p_j, \ldots, p_n)}.
\end{align*}
Let $I_{ij} \subseteq \scrO_{Z_n}$ be the ideal sheaf defining this closed
immersion, and $\Delta^{(n)}_{X/S}$ the closed subscheme of $Z_n$ defined by
$\sum_{0 \leq i,j \leq n} I_{ij}^2$.  It is easy to see that $I_{ij}$ is
generated by elements of the form
\[ \delta^{i,j} p := 1 \otimes (e^{\pi_j^* p - \pi_i^* p} - 1) \in B_n' \] for $p \in P^{\gp}$ and
\[ d^{i,j} y := (\pi_j^* y - \pi_i^* y) \otimes 1 \] for $y \in B$.  Let $\tilde
I_{ij}$ be the image of $I_{ij}$ in $\scrO_{\Delta^{(n)}_{X/S}}$.

We first note the following for future reference:

\begin{lemma}
  Let $0 \leq i, j, k, \ell \leq n$.
  \begin{enumerate}
  \item Assume $i < j$ and $k < \ell$.  Then $\mu_{ij}(I_{k\ell}) = 0$ if $i =
    k$ and $j = \ell$; otherwise, $\mu_{ij}(I_{k\ell}) = I_{k' \ell'}$, where
    \[ k' =
    \begin{cases}
      k, & k < j; \\
      i, & k = j; \\
      k - 1, & k > j,
    \end{cases} \]
    and similarly for $\ell'$.  Hence $\mu_{ij}$ gives a well-defined map
    $\scrO_{\Delta^{(n)}_{X/S}} \to \scrO_{\Delta^{(n-1)}_{X/S}}$, and the same is true with
    $\tilde I_{k\ell}$ and $\tilde I_{k' \ell'}$ in place of $I_{k\ell}$ and
    $I_{k' \ell'}$.
  \item $I_{i\ell} \subseteq I_{ij} + I_{j\ell}$.
  \end{enumerate}

  \begin{proof}
    The first statement follows from the fact that $\mu_{ij}$ acts the
    same on the generators $d^{k,\ell} y$ and $\delta^{k,\ell} p$ of $I_{k\ell}$.

    For the second statement, note that $d^{i,\ell} y = d^{i,j} y +
    d^{j,\ell} y$ for $y \in B$.  Similarly, since $1 + \delta^{i,j} p = 1 \otimes
    e^{\pi_j^* p - \pi_i^* p}$, for $p \in P^{\gp}$ we have
    \[ 1 + \delta^{i,\ell} p = (1 + \delta^{i,j} p) (1 + \delta^{j,\ell} p). \]
    Therefore, $\delta^{i,\ell} p = \delta^{i,j} p + \delta^{j,\ell} p + (\delta^{i,j} p)
    (\delta^{j,\ell} p) \in J_{ij} + J_{j\ell}$ also.
  \end{proof}
\end{lemma}

Let $\Omega^{(n)}_{X/S}$ be the $n$th antisymmetric product of
$\Omega^1_{X/S}$.  We first define a map $\bigcap_{i=1}^n \tilde I_{0i} \to
\Omega^{(n)}_{X/S}$.

\begin{proposition}
  There exists a unique $A$-linear map $\Psi_n : B_n' \to \Gamma(X,
  \Omega^{(n)}_{X/S})$ such that for $y_0, \ldots, y_n \in B$, $(p_0, \ldots, p_n) \in
  P_n'$, we have
\begin{align*}
  \Psi_n[(y_0 \otimes & \cdots \otimes y_n) \otimes e^{(p_0, \ldots, p_n)}] = \\
  & y_0 \alpha(p_0 + p_1 + \cdots + p_n) (dy_1 + y_1 \dlog
  p_1) \wedgetilde \cdots \wedgetilde (dy_n + y_n \dlog p_n).
\end{align*}
(Here $\wedgetilde$ denotes the product in the antisymmetric product
algebra $\Omega^{(\cdot)}_{X/S}$.)

\begin{proof}
  The uniqueness is clear.  To see the map is well-defined, we have
  several things to check:

  \begin{itemize}
  \item The above expression is $A$-multilinear in the variables $y_0,
    y_1, \ldots, y_n$.

    This is clear from the $A$-linearity of $d$.
  \item The expression above is independent of the choice of $p_0, \ldots,
    p_n \in P^{\gp}$.

    This follows from the fact that $p_0 + \cdots + p_n \in P$ is
    well-defined, and the fact that $\dlog$ induces a well-defined map
    $P^{\gp} / Q^{\gp} \to \Omega^1_{X/S}$.
  \item For $p' \in P$,
    \begin{align*}
      \Psi_n[(y_0 \otimes \cdots \otimes y_i \alpha(p') \otimes \cdots \otimes y_n) \otimes e^{(p_0, \ldots, p_i, \ldots, p_n)}] & =
      \\
      \Psi_n[(y_0 \otimes \cdots \otimes y_i \otimes \cdots \otimes y_n) \otimes e^{(p_0, \ldots, p_i + p', \ldots, p_n)}] &.
    \end{align*}
    For $i=0$, this is clear.  Otherwise, for $i>0$, this follows from
    the formula
    \[ d(y_i \alpha(p')) = y_i d(\alpha(p')) + \alpha(p') dy_i = \alpha(p') [dy_i + y_i
    \dlog p']. \]
  \end{itemize}
\end{proof}
\end{proposition}

\begin{remark}
  In the case of trivial log structure, i.e.~$P=0$, the
  formula for $\Psi_n$ reduces to
  \[ \Psi_n(y_0 \otimes y_1 \otimes \cdots \otimes y_n) = y_0 \, dy_1 \wedgetilde \cdots \wedgetilde
  dy_n. \]
  This is the isomorphism commonly used in synthetic differential
  geometry, for example in \cite{kock}.

  Also note that if in fact $p_0, \ldots, p_n \in P$, then
  \[ (y_0 \otimes y_1 \otimes \cdots \otimes y_n) \otimes e^{(p_0, p_1, \ldots, p_n)} =
  (y_0 \alpha(p_0) \otimes y_1 \alpha(p_1) \otimes \cdots \otimes y_n \alpha(p_n)) \otimes 1, \]
  and in this case the formula for $\Psi_n$ agrees with
  \[ y_0 \alpha(p_0) \, d(y_1 \alpha(p_1)) \wedgetilde \cdots \wedgetilde d(y_n
  \alpha(p_n)). \]
  Thus we may view the given formula for $\Psi_n$ as a natural
  generalization of the simpler formula from the case of trivial log
  structure.
\end{remark}

\begin{proposition}
  \label{prop:psinann}
  The map $\Psi_n$ annihilates $I_{ij}^2$ for each pair $0 \leq i,j \leq n$.

  \begin{proof}
    We first check the case $i=0$.  In this case, since $d(y_j y) =
    y_j \, dy + y \, dy_j$, it is straightforward to calculate that
    for $x \in B_n'$, $y \in B$, we have
    \[ \Psi_n(x \, d^{0,j} y) = (-1)^{j-1} dy \wedgetilde \Psi_{n-1}(\mu_{0,j}
    x). \]
    Therefore, if $x \in I_{0j}$, then $\Psi_n(x \, d^{0,j} y) = 0$.
    Similarly, for $p \in P^{\gp}$,
    \[ \Psi_n(x \, \delta^{0,j} p) = (-1)^{j-1} \dlog p \wedgetilde
    \Psi_{n-1}(\mu_{0,j} x), \] so again if $x \in I_{0j}$, then $\Psi_n(x \, \delta^{0,j}
    p) = 0$.

    Now for the general case, by symmetry assume $i < j$.  We observe
    that $d^{i,j} y = d^{0,j} y - d^{0,i} y$.  Thus, if $y, y' \in B$,
    then
    \begin{equation}
      \label{eq:ddprod}
      (d^{i,j} y) (d^{i,j} y') \equiv -[(d^{0,i} y) (d^{0,j} y') + (d^{0,i} y')
      (d^{0,j} y)] \pmod{J_{0i}^2 + J_{0j}^2}.
    \end{equation}
    However, since $\mu_{0j} (x(d^{0,i} y')) = (\mu_{0j} x) (d^{0,i} y)$
    for $x \in B_n'$,
    we have
    \begin{align*}
      \Psi_n(x (d^{0,i} y) (d^{0,j} y')) & = (-1)^{j-1} dy' \wedgetilde
      \Psi_{n-1}((\mu_{0j} x) (d^{0,i} y)) \\
      & = (-1)^{i+j} dy' \wedgetilde dy \wedgetilde \Psi_{n-2}(\mu_{0i} \mu_{0j}
      x).
    \end{align*}
    Therefore, for $x \in B_n'$,
    \[ \Psi_n(x (d^{i,j} y) (d^{i,j} y')) = (-1)^{i+j+1} (dy' \wedgetilde
    dy + dy \wedgetilde dy')
    \wedgetilde \Psi_{n-2}(\mu_{0i} \mu_{0j} x) = 0. \]
    Similarly, since $1 + \delta^{i,j} p = 1 \otimes e^{\pi_j^* p - \pi_i^* p}$, we have $1 +
    \delta^{0,j} p = (1 + \delta^{0,i} p) (1 + \delta^{i,j} p)$.  Multiplying both
    sides by $1 - \delta^{0,i} p$, this implies
    \[ 1 + \delta^{i,j} p \equiv (1 + \delta^{0,j} p) (1 - \delta^{0,i} p) \pmod{J_{0i}^2
      + J_{0j}^2}, \] so $\delta^{i,j} p \equiv \delta^{0,j} p - \delta^{0,i} p - (\delta^{0,j}
    p) (\delta^{0,i} p)$.  Therefore,
    \begin{equation}
      \label{eq:dlogdprod}
      (\delta^{i,j} p) (d^{i,j} y) \equiv -[(\delta^{0,i} p) (d^{0,j} y) + (d^{0,i}
      y) (\delta^{0,j} p)] \pmod{J_{0i}^2 + J_{0j}^2}
    \end{equation}
    and
    \begin{equation}
      \label{eq:dlogdlogprod}
      (\delta^{i,j} p) (\delta^{i,j} p') \equiv -[(\delta^{0,i} p) (\delta^{0,j} p') +
      (\delta^{0,i} p') (\delta^{0,j} p)] \pmod{J_{0i}^2 + J_{0j}^2}.
    \end{equation}
    From these formulas, the proof that
    $\Psi$ annihilates $x (\delta^{i,j} p) (d^{i,j} y)$ and $x (\delta^{i,j} p)
    (\delta^{i,j} p')$ proceeds as before.
  \end{proof}
\end{proposition}

Therefore, since $\Omega^{(n)}_{X/S}$ is a quasi-coherent $\scrO_X$-module,
$\Psi_n$ induces a map $\Psi_n : \scrO_{\Delta^{(n)}_{X/S}} \to \Omega^{(n)}_{X/S}$, which we
will restrict to the ideal $\bigcap_{j=1}^n \tilde I_{0j}$ of
$\scrO_{\Delta^{(n)}_{X/S}}$.  We now turn to defining a map in the other
direction.

\begin{proposition}
  \label{prop:phiiexists}
  For each $i$ with $0 < i \leq n$, there is a unique $B$-linear map $\phi_i
  : \Omega^1_{X/S} \to \tilde I_{0i}$ such that $\phi_i(dy) = d^{0,i} y$ for $y
  \in B$ and $\phi_i(\dlog p) = \delta^{0,i} p$ for $p \in P^{\gp}$.  (Here we
  consider $\tilde I_{0i}$ to be a $B$-module via $\pi_0^*$.)

  \begin{proof}
    By the universal property of $\Omega^1_{X/S}$, we need only check that
    $(D, \delta) : (B, P^{gp}) \to \tilde I_{0i}$ defined by $Dy = d^{0,i} y$
    and $\delta p = \delta^{0,i} p$ is a log derivation over $A$.  However, $D$
    is clearly $A$-linear, and since
    \[ (d^{0,i} y) (d^{0,i} y') =
    d^{0,i} (y y') - (\pi_0^* y) d^{0,i} y' - (\pi_0^* y') d^{0,i} y \in
    I_{0i}^2, \]
    $D$ is also a derivation.  Similarly, for $p \in P$, we
    have $d^{0,i} (\alpha(p)) = (\pi_0^* \alpha(p)) \delta^{0,i} p$.  Finally, to see
    that $\delta$ is additive, since $1 + \delta^{0,i} p = 1 \otimes e^{\pi_i^* p - \pi_0^* p}$,
    we have
    \[ 1 + \delta^{0,i} (p + p') = (1 + \delta^{0,i} p) (1 + \delta^{0,i}
    p') = 1 + \delta^{0,i} p + \delta^{0,i} p' + (\delta^{0,i} p) (\delta^{0,i} p'). \]
    Therefore, $\delta^{0,i} (p + p') \equiv \delta^{0,i} p + \delta^{0,i} p' \pmod{J_{0i}^2}$.
  \end{proof}
\end{proposition}

\begin{proposition}
  There is a unique map $\Phi_n : \Omega^{(n)}_{X/S} \to \prod_{j=1}^n \tilde
  J_{0j}$ such that for $\omega_1, \omega_2, \ldots, \omega_n \in \Omega^1_{X/S}$,
  \[ \Phi_n(\omega_1 \wedgetilde \omega_2 \wedgetilde \cdots \wedgetilde \omega_n) = \phi_1(\omega_1)
  \phi_2(\omega_2) \cdots \phi_n(\omega_n). \]

  \begin{proof}
    Since the formula above is clearly multilinear in $\omega_1, \ldots, \omega_n$,
    we need only check it is antisymmetric.  We claim that in fact,
    for $\omega, \tau \in \Omega^1_{X/S}$, $\phi_i(\omega) \phi_j(\tau) + \phi_i(\tau) \phi_j(\omega) = 0$ in
    $\scrO_{\Delta^{(n)}_{X/S}}$.  To see this, we refer again to the
    formulas (\ref{eq:ddprod}) through (\ref{eq:dlogdlogprod}).  Thus,
    if $\omega = dy$ and $\tau = dy'$, then by (\ref{eq:ddprod}),
    \[ (d^{0,i} y) (d^{0,j} y') + (d^{0,i} y') (d^{0,j} y) \equiv -(d^{i,j}
    y) (d^{i,j} y') \pmod{J_{0i}^2 + J_{0j}^2}, \]
    so $\phi_i(\omega) \phi_j(\tau) + \phi_i(\tau) \phi_j(\omega) \in J_{0i}^2 + J_{0j}^2 + J_{ij}^2$.
    Similarly, for the cases $\omega = \dlog p$, $\tau = dy$ and $\omega = \dlog
    p$, $\tau = \dlog p'$, we use the corresponding formulas
    (\ref{eq:dlogdprod}) and (\ref{eq:dlogdlogprod}).
  \end{proof}
\end{proposition}

We now show the two maps defined above are inverses.

\begin{theorem}
  \begin{enumerate}
  \item The composition
    \[ \Omega^{(n)}_{X/S} \overset{\Phi_n}{\longrightarrow} \prod_{j=1}^n \tilde J_{0j}
    \overset{\Psi_n}{\longrightarrow} \Omega^{(n)}_{X/S} \]
    is the identity on $\Omega^{(n)}_{X/S}$.
  \item The composition
    \[ \bigcap_{j=1}^n \tilde J_{0j} \overset{\Psi_n}{\longrightarrow} \Omega^{(n)}_{X/S}
    \overset{\Phi_n}{\longrightarrow} \prod_{j=1}^n \tilde J_{0j} \hookrightarrow \bigcap_{j=1}^n \tilde J_{0j} \]
    is the identity map on $\bigcap_{j=1}^n \tilde J_{0j}$.
  \end{enumerate}

  \begin{proof}
    For the first composition, it suffices to check for
    \[ \omega = dy_1 \wedgetilde \cdots \wedgetilde dy_i \wedgetilde
    \dlog p_{i+1} \wedgetilde \cdots \wedgetilde \dlog p_n, \]
    for $y_1, \ldots, y_i \in B$, $p_{i+1}, \ldots, p_n \in P^{\gp}$.  However,
    \begin{align*}
      \Phi_n(\omega) = \sum_{S \subseteq \{ 1, \ldots, n \}} (-1)^{|S|} & (y_{0S} \otimes y_{1S}
      \otimes \cdots \otimes y_{iS} \otimes 1 \otimes \cdots \otimes 1) \otimes \\
      & e^{(p_{0S}, 0, \ldots, 0, p_{i+1,S}, \ldots, p_{nS})},
    \end{align*}
    where:
    \begin{itemize}
    \item $y_{0S} = \prod_{1 \leq j \leq i, j \in S} y_j$;
    \item $y_{jS} = y_j$ if $j \notin S$ and $y_{jS} = 1$ if $j \in S$, $1 \leq
      j \leq i$;
    \item $p_{0S} = -\sum_{i < j \leq n, j \notin S} p_j$;
    \item $p_{jS} = p_j$ if $j \notin S$ and $p_{jS} = 0$ if $j \in S$, $i <
      j \leq n$.
    \end{itemize}
    Therefore,
    \[ \Psi_n(\Phi_n(\omega)) = \sum_S y_{0S} \, dy_{1S} \wedgetilde \cdots \wedgetilde
    dy_{iS} \wedgetilde \dlog p_{i+1,S} \wedgetilde \cdots \wedgetilde
    \dlog p_{nS}. \] However, if $S \neq \emptyset$, then either $y_{jS} = 1$ or
    $p_{jS} = 0$ for some $j \in S$, so the corresponding term is zero.
    On the other hand, for $S = \emptyset$, the corresponding term is exactly
    $\omega$.

    For the second composition, let $x = (y_0 \otimes \cdots \otimes y_n) \otimes e^{(p_0, \ldots,
      p_n)} \in B_n'$.  Then since $\pi_i^* y_i \equiv \pi_0^* y_i \pmod{J_{0i}}$,
    we calculate that
    \begin{align*}
      \phi_i(dy_i + y_i \dlog p_i) & \equiv \pi_i^* y_i \otimes 1 - \pi_0^* y_i \otimes 1 +
      \pi_i^* y_i \otimes (e^{\pi_i^* p_i - \pi_0^* p_i} - 1) \\
      & = \pi_i^* y_i \otimes
      e^{\pi_i^* p_i - \pi_0^* p_i} - \pi_0^* y_i \otimes 1 \pmod{J_{0i}^2}.
    \end{align*}
    From this we see that
    \[ \Phi_n \circ \Psi_n = \prod_{j=1}^n (\id - M_j) \]
    on $\scrO_{\Delta^{(n)}_{X/S}}$, where
    \begin{align*}
      M_j[ & (y_0 \otimes y_1 \otimes \cdots \otimes y_j \otimes \cdots \otimes y_n) \otimes e^{(p_0, p_1, \ldots, p_j, \ldots, p_n)}] = \\
      & (y_0 y_j \otimes y_1 \otimes \cdots \otimes 1 \otimes \cdots \otimes y_n) \otimes e^{(p_0 + p_j, p_1, \ldots, 0, \ldots, p_n)}.
    \end{align*}
    However, $M_j$ factors through $\mu_{0j}$, so this implies that $\Phi_n
    \circ \Psi_n = \id$ on $\bigcap_{j=1}^n \tilde I_{0j}$.
  \end{proof}
\end{theorem}

\begin{corollary}
  We have
  \[ \bigcap_{j=1}^n \tilde J_{0j} = \prod_{j=1}^n \tilde J_{0j} = \bigcap_{0 \leq i,j \leq
    n} \tilde J_{ij} \simeq \Omega^{(n)}_{X/S}. \]

  \begin{proof}
    From the theorem, the inclusion map $\prod_{j=1}^n \tilde J_{0j} \hookrightarrow
    \bigcap_{j=1}^n \tilde J_{0j}$ must be surjective, so it is the identity
    map and the two ideals are equal.  Thus $\Phi_n$ and $\Psi_n$ are
    inverse isomorphisms between this common ideal and
    $\Omega^{(n)}_{X/S}$.  Now clearly, $\bigcap_{0 \leq i,j \leq n} \tilde J_{ij} \subseteq
    \bigcap_{j=1}^n \tilde J_{0j}$.  On the other hand, if $0 \leq i, j \leq n$,
    then $\tilde J_{0j} \subseteq \tilde J_{0i} + \tilde J_{ij}$, so $\tilde
    J_{0i} \tilde J_{0j} \subseteq \tilde J_{ij}$.  Therefore, $\prod_{j=1}^n
    \tilde J_{0j} \subseteq \bigcap_{0 \leq i,j \leq n} \tilde J_{ij}$ also.
  \end{proof}
\end{corollary}

We now present another formulation of this ideal which is more useful
in certain situations.  First, we note that
\begin{align*}
  \tilde J_{01} & = \tilde J_{01}; \\
  \tilde J_{02} & \subseteq \tilde J_{01} + \tilde J_{12}; \\
  \tilde J_{03} & \subseteq \tilde J_{01} + \tilde J_{12} + \tilde J_{23}; \\
  & \vdots \\
  \tilde J_{0n} & \subseteq \tilde J_{01} + \tilde J_{12} + \cdots + \tilde J_{n-1,n}.
\end{align*}
Therefore, $\prod_{j=1}^n \tilde J_{0j} \subseteq \prod_{j=1}^n \tilde J_{j-1,j}$.
Similarly, since $\tilde J_{j-1,j} \subseteq \tilde J_{0,j-1} + \tilde
J_{0,j}$, the reverse inclusion also holds, and $\prod_{j=1}^n \tilde
J_{0j} = \prod_{j=1}^n \tilde J_{j-1,j}$.

In fact, since $d^{0,j} y \equiv d^{j-1,j} y \pmod{J_{0,j-1}}$ and $\delta^{0,j}
p \equiv \delta^{j-1,j} p \pmod{J_{0,j-1}}$, while $J_{0,j-1} \subseteq J_{01} + J_{12} +
\cdots + J_{j-2,j-1}$, we see that
\[ \Phi_n(\omega_1 \wedgetilde \cdots \wedgetilde \omega_n) = \psi_1(\omega_1) \cdots \psi_n(\omega_n), \]
where $\psi_i(dy) = d^{i-1,i} y$ and $\psi_i(\dlog p) = \delta^{i-1,i} p$.  From
this we may calculate that
\[ \Phi_n \circ \Psi_n = (\id - M_n') \circ \cdots \circ (\id - M_2') \circ (\id - M_1'). \]
Here
\begin{align*}
  M_j'[ & (y_0 \otimes \cdots \otimes y_{j-1} \otimes y_j \otimes \cdots \otimes y_n) \otimes e^{(p_0, \ldots, p_{j-1}, p_j,
  \ldots, p_n)}] \\
  = & (y_0 \otimes \cdots \otimes y_{j-1} y_j \otimes 1 \otimes \cdots \otimes y_n) \otimes e^{(p_0, \ldots, p_{j-1} + p_j,
  0, \ldots, p_n)}.
\end{align*}
However, $M_j'$ factors through $\mu_{j-1,j}$, hence $\Phi_n \circ \Psi_n = \id$
on $\bigcap_{j=1}^n \tilde J_{j-1,j}$.  From this we conclude that
\[ \bigcap_{j=1}^n \tilde J_{j-1,j} = \prod_{j=1}^n \tilde J_{j-1,j} = \bigcap_{j=1}^n
\tilde J_{0j} = \prod_{j=1}^n \tilde J_{0j} = \bigcap_{0 \leq i,j \leq n} \tilde
J_{ij}. \]

\subsection{Globalization}
Now let $X \to S$ be an arbitrary morphism of fine log schemes, and let
$\Delta^n_{X/S} := X \times_S X \times_S \cdots \times_S X$ be the $n+1$-fold product.  Let $\tilde
\Delta^n_{X/S}$ be the log formal neighborhood of the diagonal immersion $\Delta : X \to
\Delta^n_{X/S}$.  Then for $0 \leq i < j \leq n$, we have a closed immersion $m_{ij} :
\Delta^{n-1}_{X/S} \to \Delta^n_{X/S}$ defined by
\begin{align*}
  m_{ij} & (x_0, \ldots, x_i, \ldots, x_{j-1}, x_j, \ldots, x_{n-1}) = \\
  & (x_0, \ldots, x_i, \ldots, x_{j-1}, x_i, x_j, \ldots, x_{n-1}).
\end{align*}
This induces a closed immersion $\tilde \Delta^{n-1}_{X/S} \to \tilde \Delta^n_{X/S}$.  Let
$J_{ij}$ be the ideal of $\scrO_{\tilde \Delta^n_{X/S}}$ defining this closed
immersion, and let $\Delta^{(n)}_{X/S}$ be the closed subscheme of $\tilde \Delta^n_{X/S}$
defined by $\sum_{0 \leq i < j \leq n} J_{ij}^2$.  Finally, let $\tilde J_{ij}$
be the image of $J_{ij}$ in $\scrO_{\Delta^{(n)}_{X/S}}$.

\begin{theorem}
  We have
  \[ \bigcap_{j=1}^n \tilde J_{j-1,j} = \prod_{j=1}^n \tilde J_{j-1,j} =
  \bigcap_{j=1}^n \tilde J_{0j} = \prod_{j=1}^n \tilde J_{0j} = \bigcap_{0 \leq i,j \leq n}
  \tilde J_{ij}, \]
  and this ideal is canonically isomorphic to $\Omega^{(n)}_{X/S}$.

  \begin{proof}
    Since the construction above is local with respect to both $X$ and
    $S$, in proving the first statement we may assume that $X$ and $S$
    are affine and that we have a chart $(P \to \scrO_X, Q \to \scrO_S, Q
    \to P)$ of the morphism $X \to S$.  Then note that in the construction
    of the previous section, the ideal $J$ of the diagonal immersion
    $X \to \Delta^n_{X/S}$ satisfies
    \[ J \subseteq J_{01} + J_{12} + \cdots + J_{n-1,n}. \]
    Therefore, $J^{n+1} \subseteq \sum_{0 \leq i,j \leq n} J_{ij}^2$, and in fact we
    could have used the $n$th log infinitesimal neighborhood of the
    diagonal in place of $Z_n$.  The same holds true for the global
    construction above.  Now it is easy to see that the global
    construction reduces to the local construction of the last section
    in this case.  From this we immediately see the equality of the
    five ideals.

    Now to establish an isomorphism between $\prod_j \tilde J_{0j}$ and
    $\Omega^{(n)}_{X/S}$, a similar proof to the proof of
    \ref{prop:phiiexists} shows that there are unique maps $\phi_i :
    \Omega^1_{X/S} \to \tilde J_{0i}$ such that $\phi_i(dy) = \pi_i^* y - \pi_0^* y$
    for $y \in \scrO_X$, and $\phi_i(\dlog m) = \alpha(\pi_i^* m - \pi_0^* m) - 1$
    for $m \in \MXgp$.  (Here, since $X \to \tilde \Delta^n_{X/S}$ is exact and
    $\pi_i^* m - \pi_0^* m$ pulls back to 0 in $\scrM_X^{\gp}$, we must
    have $\pi_i^* m - \pi_0^* m \in \scrM_{\tilde \Delta^n_{X/S}}$.)  Therefore,
    there is a unique map $\Phi_n : \Omega^{(n)}_{X/S} \to \prod_{j=1}^n \tilde
    J_{0j}$ such that
    \[ \Phi_n(\omega_1 \wedgetilde \cdots \wedgetilde \omega_n) = \phi_1(\omega_1) \cdots \phi_n(\omega_n) \]
    for $\omega_1, \ldots, \omega_n \in \Omega^1_{X/S}$.  (The map exists locally by the
    previous section, and the uniqueness allows us to glue the local
    maps.)  By the previous section, $\Phi_n$ is locally an isomorphism,
    so $\Phi_n$ gives a global isomorphism $\Omega^{(n)}_{X/S} \overset{\sim}{\to}
    \prod_{j=1}^n \tilde J_{0j}$.
  \end{proof}
\end{theorem}

\subsection{The Divided Power Envelope}

In this section, let $D(n)$ denote the log PD envelope of the diagonal
in $\Delta^n_{X/S}$.  As before we get closed immersions $m_{ij} : D(n-1) \to
D(n)$.  Let $\bar J_{ij} \subseteq \scrO_{D(n)}$ be the PD ideal corresponding
to $m_{ij}$, $\Delta^{[n]}_{X/S}$ be the closed subscheme of $D(n)$ defined
by $\sum_{i,j} \bar J_{ij}^{[2]}$, and $\tilde J_{ij}$ be the ideal of
$\scrO_{\Delta^{[n]}_{X/S}}$ corresponding to $\bar J_{ij}$.

\begin{theorem}
  \[ \bigcap_{j=1}^n \tilde J_{0j} = \prod_{j=1}^n \tilde J_{0j} = \bigcap_{j=1}^n
  \tilde J_{j-1,j} = \prod_{j=1}^n \tilde J_{j-1,j} = \bigcap_{0 \leq i,j \leq n}
  \tilde J_{ij}, \]
  and this ideal is canonically isomorphic to $\Omega^n_{X/S}$.

  \begin{proof}
    First, observe that by the universal property of the PD envelope,
    $D(n) \simeq D(1) \times_X \cdots \times_X D(1)$, the product of $n$ factors of
    $D(1)$, each considered as a scheme over $X$ via the projection
    $\pi_0$.  Therefore, $\Delta^{[1]}_{X/S} \times_X \cdots \times_X \Delta^{[1]}_{X/S}$ is isomorphic to
    the closed subscheme of $D(n)$ corresponding to $\sum_j \bar
    J_{0j}^{[2]}$.  Also, it is easy to see that $\Delta^{[1]}_{X/S} \simeq \Delta^{(1)}_{X/S}$.

    Therefore, to see that the previous map $\Psi_n$ induces a map
    $\scrO_{\Delta^{[n]}_{X/S}} \to \Omega^n_{X/S}$, it suffices to check that for $1 \leq
    i,j \leq n$, $\Psi_n[x \cdot (d^{i,j} y)^{[2]}] = \Psi_n[x \cdot (\delta^{i,j} p)^{[2]}]
    = 0$ for $x \in \scrO_{D(n)}$, $y \in \scrO_X$, $p \in \MXgp$.  However,
    since $d^{i,j} y = d^{0,j} y - d^{0,i} y$, we have
    \begin{align*}
      (d^{i,j} y)^{[2]} & = (d^{0,j} y)^{[2]} + (d^{0,i} y)^{[2]} -
      (d^{0,i} y) (d^{0,j} y) \\
      & \equiv -(d^{0,i} y) (d^{0,j} y) \pmod{\bar J_{0i}^{[2]} + \bar
        J_{0j}^{[2]}}.
    \end{align*}
    Therefore,
    \begin{align*}
      \Psi_n(x \cdot (d^{i,j} y)^{[2]}) & = -\Psi_n(x (d^{0,i} y) (d^{0,j} y))
      \\
      & = (-1)^{i+j} dy \land dy \land \Psi_{n-2}(\mu_{0i} \mu_{0j} x) = 0.
    \end{align*}
    Similarly,
    \[ (\delta^{i,j} p)^{[2]} \equiv -(\delta^{0,i} p) (\delta^{0,j} p) \pmod{\bar J_{0i}^{[2]} +
      \bar J_{0j}^{[2]}}. \]
    Therefore,
    \begin{align*}
      \Psi_n(x \cdot (\delta^{i,j} p)^{[2]}) & = -\Psi_n(x (\delta^{0,i} p) (\delta^{0,j} p))
      \\
      & = (-1)^{i+j} \dlog p \land \dlog p \land \Psi_{n-2} (\mu_{0i} \mu_{0j} x) =
      0.
    \end{align*}

    Similarly, we already know $\Phi_n : \Omega^{(n)}_{X/S} \to \prod_{j=1}^n \tilde
    J_{0j}$ is antisymmetric.  Therefore, to check $\Phi_n$ induces a map
    $\Omega^n_{X/S} \to \prod_{j=1}^n \tilde J_{0j}$, it suffices to check that
    it annihilates $dy \wedgetilde dy \wedgetilde \omega$ and $\dlog p
    \wedgetilde \dlog p \wedgetilde \omega$.  However, from the above, we
    see that in fact \[ (d^{0,1} y) (d^{0,2} y) \equiv -(d^{1,2} y)^{[2]}
    \pmod{\bar J_{01}^{[2]} + \bar J_{02}^{[2]}}, \] so $\phi_1(dy)
    \phi_2(dy) \in \bar J_{01}^{[2]} + \bar J_{02}^{[2]} + \bar
    J_{12}^{[2]}$, and this is zero in $\scrO_{\Delta^{[n]}_{X/S}}$.
    The proof that $\phi_1(\dlog p) \phi_2(\dlog p) = 0$ in
    $\scrO_{\Delta^{[n]}_{X/S}}$ is similar.

    Now since $\Psi_n$ and $\Phi_n$ were induced from inverse maps, they are
    inverse isomorphisms, and the equality of the ideals follows as
    before.
  \end{proof}
\end{theorem}

%% file: derham.tex
\section{The de~Rham Complex}

We now describe the de~Rham complex $\Omega^\cdot_{X/S}$ in terms of our
characterization of $\Omega^n_{X/S}$.  We begin with the wedge product:
thus, let $m, n > 0$.  Then we have a map
\begin{align*}
  \Delta^{m+n}_{X/S} & \to \Delta^m_{X/S} \times_X \Delta^n_{X/S},\\
  (x_0, \ldots, x_m, \ldots, x_{m+n}) & \mapsto ((x_0, \ldots, x_m), (x_m, \ldots, x_{m+n})).
\end{align*}
Here we consider $\Delta^m_{X/S}$ as a scheme over $X$ via the last projection
$\pi_m$, and $\Delta^n_{X/S}$ as a scheme over $X$ via the first projection $\pi_0$.
This induces a map $D(m+n) \to D(m) \times_X D(n)$, which in turn induces a
map
\[ s_{mn} : \Delta^{[m+n]}_{X/S} \to \Delta^{[m]}_{X/S} \times_X \Delta^{[n]}_{X/S}. \]
Locally, this map is also induced by the ``smashing'' map
\begin{align*}
  B_m' \otimes_B B_n' \to & B_{m+n}', \\
  {} [(y_0 \otimes \cdots \otimes y_m) \otimes e^{(p_0, \ldots, p_m)}] \otimes [(y_0' & \otimes \cdots \otimes y_n') \otimes
  e^{(p_0', \ldots, p_n')}] \mapsto \\
  (y_0 \otimes \cdots \otimes y_m y_0' \otimes \cdots \otimes y_n') & \otimes e^{(p_0, \ldots, p_m + p_0', \ldots, p_n')}.
\end{align*}

\begin{remark}
  Although we also have a map $D(m) \times_X D(n) \to D(m+n)$, we do not get
  an induced map $\Delta^{[m]}_{X/S} \times_X \Delta^{[n]}_{X/S} \to \Delta^{[m+n]}_{X/S}$
  in general.  For example, the pullback of $\bar J_{0,m+n}^{[2]}$
  does not correspond to anything from $\sum_{0\leq i,j\leq m} \bar J_{ij}^{[2]}$
  or $\sum_{0\leq i,j\leq n} \bar J_{ij}^{[2]}$.
\end{remark}

\begin{proposition}
  We have a commutative diagram
  \[
  \begin{CD}
    \bigcap_{i=1}^m \tilde J_{i-1,i} \otimes_{\scrO_X} \bigcap_{j=1}^n \tilde J_{j-1,j}
    @> s_{mn}^* >> \bigcap_{j=1}^{m+n} \tilde J_{j-1,j} \\
    @V \Psi_m \otimes \Psi_n V \simeq V @V \Psi_{m+n} V \simeq V \\
    \Omega^m_{X/S} \otimes_{\scrO_X} \Omega^n_{X/S} @> \land >> \Omega^{m+n}_{X/S}.
  \end{CD}
  \]

  \begin{proof}
    Our first task is to verify that $s_{mn}^*$ actually induces a map
    as in the top row.  To see this, note that
    \[ s_{mn} \circ m_{j-1,j} =
    \begin{cases}
      (m_{j-1,j}, \id) \circ s_{m-1,n}, & j \leq m; \\
      (\id, m_{j-m-1, j-m}) \circ s_{m,n-1}, & j > m.
    \end{cases}
    \]
    Therefore, converting to dual statements in terms of $s_{mn}^*$
    and $m_{j-1,j}^*$, we see that the image of $s_{mn}^*$ is
    annihilated by each $m_{j-1,j}^*$ and is thus in each kernel
    $\tilde J_{j-1,j}$.

    Now to check the commutativity, we first reverse the vertical
    arrows and replace them by $\Phi_m \otimes \Phi_n$ and $\Phi_{m+n}$,
    respectively.  Now, from the fact that
    \[ \Phi_n(\omega_1 \land \cdots \land \omega_n) = \psi_1(\omega_1) \cdots \psi_n(\omega_n), \]
    where $\psi_j(dy) = d^{j-1,j} y$ and $\psi_j(\dlog p) = \delta^{j-1,j} p$,
    the commutativity is clear.
  \end{proof}
\end{proposition}

We now turn to the differential map in the de Rham complex; thus, fix
$n \geq 0$.  Then for $0 \leq j \leq n+1$, we have maps
\[ \Delta^{n+1}_{X/S} \to \Delta^n_{X/S}, ~ (x_0, \ldots, x_j, \ldots, x_{n+1}) \mapsto (x_0, \ldots, \hat x_j, \ldots,
x_{n+1}). \]
These induce maps $D(n+1) \to D(n)$, which in turn induce maps
\[ d_j : \Delta^{[n+1]}_{X/S} \to \Delta^{[n]}_{X/S}. \]
Locally, these maps are also
induced by the insertion maps
\begin{align*}
  B_n' & \to B_{n+1}', \\
  (y_0 \otimes \cdots \otimes y_n) \otimes e^{(p_0, \ldots, p_n)} & \mapsto (y_0 \otimes \cdots \otimes 1 \otimes \cdots \otimes y_n) \otimes
  e^{(p_0, \ldots, 0, \ldots, p_n)},
\end{align*}
with insertion in the $j$th position.

\begin{proposition}
  We have a commutative diagram
  \[
  \begin{CD}
    \bigcap_{j=1}^n \tilde J_{j-1,j} @> d_0^* - d_1^* + \cdots + (-1)^{n+1} d_{n+1}^* >>
    \bigcap_{j=1}^{n+1} \tilde J_{j-1,j} \\
    @V \Psi_n V \simeq V @V \Psi_{n+1} V \simeq V \\
    \Omega^n_{X/S} @> d >> \Omega^{n+1}_{X/S}.
  \end{CD}
  \]

  \begin{proof}
    Let $e_n := d_0^* - d_1^* + \cdots + (-1)^{n+1} d_{n+1}^* :
    \scrO_{\Delta^{[n]}_{X/S}} \to \scrO_{\Delta^{[n+1]}_{X/S}}$.  Again, we first
    need to check that $e_n$ induces a map as in the top
    row.  To see this, note that
    \[ d_j \circ m_{i-1,i} =
    \begin{cases}
      m_{i-2,i-1} \circ d_j, & j < i-1; \\
      \id, & j = i-1 \text{ or } i; \\
      m_{i-1,i} \circ d_{j-1}, & j > i.
    \end{cases}
    \]
    From this, it is easy to check that the image of $e_n$ is
    annihilated by each $m_{i-1,i}^*$.

    Now it follows formally from the appropriate identities that
    $e_{n+1} \circ e_n = 0$, corresponding to the requirement that $d \circ d
    = 0$.  Furthermore,
    \[ e_{m+n} \circ s_{mn}^* = s_{m+1,n}^* \circ (e_m \otimes \id) + (-1)^m
    s_{m,n+1}^* \circ (\id \otimes e_n), \] which corresponds to the requirement
    that $d(\omega \land \tau) = d\omega \land \tau + (-1)^m \omega \land d\tau$.  It is easy to see that
    $e_0$ agrees with $d : \scrO_X \to \Omega^1_{X/S}$.  Therefore, all that
    is left is to verify that $e_1(\dlog p) = 0$ for $p \in \MXgp$.  We
    calculate locally, where $\dlog p = 1 \otimes (e^{(-p, p)} - 1) \in B_1'$ for
    $p \in P^{\gp}$.  Thus,
    \[ e_1(\dlog p) = 1 \otimes e^{(0, -p, p)} - 1 \otimes e^{(-p, 0, p)} +
    1 \otimes e^{(-p, p, 0)} - 1 \otimes 1. \]
    Now by the definition of $\Psi_2$, the last three terms are
    annihilated, and the first gets mapped to $-\dlog p \land \dlog p =
    0$.
  \end{proof}
\end{proposition}

\begin{remark}
  In the case of trivial log structure, we have an easier proof:
  we see that $\Psi_{n+1}$ annihilates the image of $d_j^*$ for $j > 0$,
  while locally,
  \[ \Psi_{n+1}[d_0^* (y_0 \otimes \cdots \otimes y_n)] = dy_0 \land \cdots \land dy_n = d [\Psi_n(y_0 \otimes \cdots
  \otimes y_n)]. \]
  However, to extend this proof to the case of log schemes, we must
  verify that
  \begin{align*}
    d[y_0 \alpha(p_0 + \cdots + p_n) (dy_1 + y_1 \dlog p_1) \land \cdots \land (dy_n + y_n
    \dlog p_n) = \\
    \alpha(p_0 + \cdots + p_n) (dy_0 + y_0 \dlog p_0) \land \cdots \land (dy_n +
    y_n \dlog p_n).
  \end{align*}
  While this can be done, we prefer to give the more conceptual proof
  above.
\end{remark}

\begin{remark}
  By taking the corresponding maps on the antisymmetric powers
  $\Omega^{(\cdot)}_{X/S}$ of $\Omega^1_{X/S}$, we can define a natural complex.
  However, from the above calculations, we see that we get $d(\dlog m)
  = \dlog m \wedgetilde \dlog m$, instead of 0.  This illustrates why
  in defining the logarithmic de~Rham complex such that $d(\dlog m) =
  0$, we need the full alternating product instead of just the
  antisymmetric product.  (This requirement appears in the need to
  check that $d^2 \alpha(m) = d(\alpha(m) \dlog m) = 0$.)
\end{remark}


%% file: coeff.tex
\section{Coefficients}

In this section, we will assume that $X$ is log smooth over $S$, and
$(E, \nabla)$ is an $\scrO_X$-module with quasi-nilpotent integrable
connection.  Then this corresponds to a crystal $E$ of
$\scrO_{X/S}$-modules on the log crystalline site $(X/S)_{cris}$.
Recall that an object of $(X/S)_{cris}$ is a tuple $(U, T, i, \delta)$
where $U$ is an open subscheme of $X$, $i : U\to T$ is an {\em exact}
log closed immersion, and $\delta$ is a PD structure on the ideal of $i$.
(For convenience of notation, we often use $T$ to represent this
object.)  Then a morphism $g : T_1 \to T_2$ in this site is a morphism
respecting the closed immersions and the PD structures, and a covering
$(U_\lambda, T_\lambda, i_\lambda, \delta_\lambda)_{\lambda \in \Lambda}$ of $T$ is a family such that $(T_\lambda)$ is
a Zariski open covering of $T$.  Giving a sheaf $E$ on this site is
then equivalent to giving a sheaf $E_T$ on $T$ for each object $T$ of
$(X/S)_{cris}$, along with transition maps $\theta_g : g^{-1} E_{T_2} \to
E_{T_1}$ for each morphism $g : T_1 \to T_2$ in the site, satisfying the
compatibility relation
\[ \theta_{hg} = \theta_g \circ g^{-1} \theta_h \] for the composition of $g : T_1 \to
T_2$, $h : T_2 \to T_3$.  We define $\scrO_{X/S}$ to be the sheaf with
$(\scrO_{X/S})_T := \scrO_T$, and a sheaf $E$ of $\scrO_{X/S}$-modules
is a {\em crystal} if for each morphism $g : T_1 \to T_2$ in
$(X/S)_{cris}$, the induced transition map $\theta_g : g^* E_{T_2} \to
E_{T_1}$ is an isomorphism.  For more details, see \cite{kato}.

We note that by construction, each $\Delta^{[n]}_{X/S}$ is an object of
$(X/S)_{cris}$, and each $d_j : \Delta^{[n+1]}_{X/S} \to \Delta^{[n]}_{X/S}$ is a
morphism in this site.  We thus get transition maps
\[ \theta_{d_j} : d_j^* E_{\Delta^{[n]}_{X/S}} \overset{\sim}{\to} E_{\Delta^{[n+1]}_{X/S}}. \]
Here we will consider $E_{\Delta^{[n]}_{X/S}}$ as being identified with $E
\otimes_{\scrO_X} \scrO_{\Delta^{[n]}_{X/S}}$ via $\theta_{\pi_0} : \pi_0^* E_X \overset{\sim}{\to}
E_{\Delta^{[n]}_{X/S}}$.

Also, since the map \[ inc \circ \Phi_n : \Omega^n_{X/S} \overset{\sim}{\to} \bigcap_{j=1}^n
\tilde J_{j-1,j} \hookrightarrow \scrO_{\Delta^{[n]}_{X/S}} \] is a split injection (with
splitting $\Psi_n$), so is the map $\id \otimes (inc \circ \Phi_n) : E \otimes_{\scrO_X}
\Omega^n_{X/S} \to E \otimes_{\scrO_X} \scrO_{\Delta^{[n]}_{X/S}}$.  Furthermore, we see that
the image is equal to the intersection of the kernels of the
transition maps $\theta_{m_{j-1,j}} : E_{\Delta^{[n]}_{X/S}} \to E_{\Delta^{[n+1]}_{X/S}}$.  We
will treat $\id \otimes \Phi_n$ as identifying $E \otimes_{\scrO_X} \Omega^n_{X/S}$ with
this submodule of $E \otimes_{\scrO_X} \scrO_{\Delta^{[n]}_{X/S}}$.

We now give a characterization of the de~Rham complex with
coefficients in $E$, in terms of the transition maps $\theta_{d_j}$ and the
above identifications.

\begin{proposition}
  We have a commutative diagram
  \[
  \begin{CD}
    @> \nabla >> E \otimes_{\scrO_X} \Omega^n_{X/S} @> \nabla >> E \otimes_{\scrO_X}
    \Omega^{n+1}_{X/S} @> \nabla >> \\
    @. @VVV @VVV @. \\
    @>>> E \otimes_{\scrO_X} \scrO_{\Delta^{[n]}_{X/S}} @>>> E \otimes_{\scrO_X}
    \scrO_{\Delta^{[n+1]}_{X/S}} @>>> \\
    @. @V \theta_{\pi_0} V \simeq V @V \theta_{\pi_0} V \simeq V @. \\
    @>>> E_{\Delta^{[n]}_{X/S}} @> \theta_{d_0} - \theta_{d_1} + \cdots +
    (-1)^{n+1} \theta_{d_{n+1}} >> E_{\Delta^{[n+1]}_{X/S}} @>>> .
  \end{CD}
  \]

  \begin{proof}
    Let $e \in E$, $\omega \in \Omega^n_{X/S} \subseteq \scrO_{\Delta^{[n]}_{X/S}}$.  If $j > 0$,
    then since the composition of $d_j : \Delta^{[n+1]}_{X/S} \to
    \Delta^{[n]}_{X/S}$ and $\pi_0 : \Delta^{[n]}_{X/S} \to X$ is equal to $\pi_0 :
    \Delta^{[n+1]}_{X/S} \to X$, we see that $\theta_{d_j}(e \otimes \omega) = e \otimes d_j^* \omega$.
    Therefore, all we need to do to finish the proof is to show that
    $\theta_{d_0}(e \otimes \omega) = e \otimes d_0^* \omega + \theta_{s_{1n}}(\nabla e \otimes \omega)$.  Then since
    $\pi_0 \circ s_{1n} = \pi_0$, $\theta_{s_{1n}}(\nabla e \otimes \omega) = \nabla e \land \omega$, and it
    immediately follows that for $\omega \in \Omega^n_{X/S}$,
    \[ (\theta_{d_0} - \theta_{d_1} + \cdots \pm \theta_{d_{n+1}}) (e \otimes \omega) = \nabla e \land \omega + e \otimes
    d\omega = \nabla(e \otimes \omega). \]

    However, by definition,
    \[ \nabla = \theta_{\pi_1} - \theta_{\pi_0} : E \to E
    \otimes_{\scrO_X} \Omega^1_{X/S} \hookrightarrow E \otimes_{\scrO_X} \scrO_{\Delta^{[1]}_{X/S}}. \]
    Hence by the linearity of $\theta_{\pi_1}$, $\theta_{\pi_1}(e \otimes \omega) = e \otimes
    \omega + (\nabla e) \omega$ for $e \in E, \omega \in \scrO_{\Delta^{[1]}_{X/S}}$.  Now considering
    the map $(\pi_0, \pi_1) : \Delta^{[n]}_{X/S} \to \Delta^{[1]}_{X/S}$, we must have $\theta_{\pi_1} =
    \theta_{(\pi_0, \pi_1)} \circ (\pi_0, \pi_1)^* \theta_{\pi_1}$, so
    \[ \theta_{\pi_1}(e \otimes \omega) = e \otimes \omega + [\theta_{(\pi_0, \pi_1)} (\nabla e)] \omega \]
    for $e \in E, \omega \in \scrO_{\Delta^{[n]}_{X/S}}$.  But since $\pi_0 \circ d_0 = \pi_1$, we
    now get
    \begin{align*}
      \theta_{d_0}(e \otimes \omega) & = \theta_{\pi_1}(e \otimes d_0^* \omega) = e \otimes d_0^* \omega +
      [\theta_{(\pi_0, \pi_1)} (\nabla e)] d_0^* \omega \\
      & = e \otimes d_0^* \omega + \theta_{s_{1n}}(\nabla e \otimes \omega).
    \end{align*}
  \end{proof}
\end{proposition}

\begin{remark}
  It is easy to see, independently of the above calculation, that $\sum_j
  (-1)^j \theta_{d_j}$ induces maps $E \otimes_{\scrO_X} \Omega^n_{X/S} \to E
  \otimes_{\scrO_X} \Omega^{n+1}_{X/S}$ forming the differential maps in a
  complex.  We thus have an alternate proof that the standard formula
  $\nabla(e \otimes \omega) = \nabla e \land \omega + e \otimes d\omega$ gives a well-defined complex, and in
  fact we see in this way that this is a natural generalization of the
  usual de Rham complex $(\Omega^\cdot_{X/S}, d)$.
\end{remark}